\documentclass[12pt,reqno]{amsart}
\usepackage{a4wide}

\numberwithin{equation}{section}

\newtheorem{theorem}{Theorem}[section]
\newtheorem{proposition}[theorem]{Proposition}

\newtheorem{lemma}[theorem]{Lemma}

\theoremstyle{definition}

\newtheorem{remark}[theorem]{Remark}

\newcommand{\R}{\mathbb{R}}

\begin{document}

\title[ Infinitely many positive solutions
 ]
{ Infinitely many positive solutions for the nonlinear Schr\"odinger equations
in $ \R^N$
 }

 \author{Juncheng Wei  and Shusen Yan
}

\address{Department of Mathematics, Chinese University of Hong Kong Shatin, Hong Kong}

\email{wei@math.cuhk.edu.hk}

\address{ School of Mathematics,Statistics and Computer
Science, The  University of New England Armidale, NSW 2351,
Australia}

\email{ syan@turing.une.edu.au }

\begin{abstract}

We  consider the following  nonlinear problem in $\R^N$
\begin{equation}
\label{eq}
- \Delta u +V(|y|)u=u^{p},\quad u>0 \ \   \mbox{in} \  \R^N, \ \ \
 u \in H^1(\R^N)
\end{equation}
where $V(r)$ is a  positive function, $1<p <\frac{N+2}{N-2}$. We
show that if $V(r)$ has the following expansion:  There  are constants $a>0$, $m>1$, $\theta>0$, and $V_0>0$,
such that
\[
V(r)= V_0+\frac a {r^m} +O\bigl(\frac1{r^{m+\theta}}\bigr),\quad \text{as $r\to +\infty$,}
\]
then \eqref{eq} has {\bf infinitely many  non-radial positive } solutions,
whose energy can be made arbitrarily large.

\end{abstract}

\maketitle

\section{Introduction}

Standing waves for the following nonlinear Schr\"odinger equation in
$\R^N$:
\begin{align}\label{nls}
-i \frac{\partial \psi}{\partial t}=\Delta \psi - \tilde V(y) \psi +
|\psi|^{p-1} \psi,
\end{align}
where $p > 1$, are solutions of the form $\psi(t,y) = \exp(i \lambda
t)u(y)$. Assuming that the amplitude $u(y)$ is positive and vanishes
at infinity, we see that $\psi$ satisfies \eqref{nls} if and only if
$u$ solves the nonlinear elliptic problem
\begin{align}\label{equation}
-\Delta u + V(y) u = u^p, \quad u>0, \quad \lim_{|y|\to+\infty}
u(y)=0,
\end{align}
where $V(y) = \tilde V(y) + \lambda$. In the rest of this paper,  we
will assume that $V$ is bounded, and $V(y)\ge V_0>0$.

A problem which is similar to \eqref{equation} is the following
scalar field equation:
\begin{align}\label{equation1}
-\Delta u +  u = Q(y)u^p, \quad u>0, \quad \lim_{|y|\to+\infty}
u(y)=0,
\end{align}
where $Q(y)$ is bounded, and $Q(y)\ge Q_0>0$.

If
\begin{equation}\label{1-18-4}
\inf_{y \in \R^N} V(y) < \lim_{|x|\to\infty} V(x),  \quad(\text{or}\;  \sup_{y \in \R^N}  Q(y)>
\lim_{|x|\to\infty} Q(x) ),
\end{equation}
then, using the concentration compactness principle \cite{L1,L2},
one can show that \eqref{equation} and \eqref{equation1} have a
least energy solution. See for example \cite{DN,L1,L2,R}. But if
\eqref{1-18-4} does not hold, \eqref{equation} (or
\eqref{equation1}) may not have least energy solution. So, one needs
to find solution with higher energy level. For results on this
aspect, the readers can refer to \cite{BL,BLi,C}. Note that the
energy of the solutions in \cite{BL,C} is less than twice of the
first level at which the  Palais-Smale condition fails.

Recently, Cerami, Devillanova and  Solimini \cite{CDS} showed that
the following problem
\[
-\Delta u + V(y) u = |u|^{p-1}u,  \quad \lim_{|y|\to+\infty}
u(y)=0,
\]
has {\bf infinitely many sign-changing solutions} if $V(y)$ tends to its
limit at infinity from below with a suitably rate. Except \cite{CR},
where $Q(y)$ is periodic, there is no result on the multiplicity of
positive solutions for \eqref{equation} (or \eqref{equation1}).

On the other hand, if we consider the following singularly perturbed
problem:
\begin{equation}\label{2-18-4}
-\varepsilon^2\Delta u + V(y) u = u^p, \quad u>0, \quad \lim_{|y|\to+\infty}
u(y)=0,
\end{equation}
or
\begin{equation}\label{3-18-4}
-\varepsilon^2\Delta u + u = Q(y)u^p, \quad u>0, \quad \lim_{|y|\to+\infty}
u(y)=0,
\end{equation}
where $\varepsilon>0$ is a small parameter, then the number of the
critical points of $V(y)$ (or $Q(y)$) (see for example
\cite{ABC,CNY1,CNY2},\cite{DF1}--\cite{DF4},\cite{FW,O,W}), the type
of the critical points of $V(y)$ (or $Q(y)$) (see for example
\cite{DY,KW,NY}, and the topology of the level set of $V(y)$ (or
$Q(y)$) \cite{AMN1, AMN2,DLY, FM}, can affect the number of the
solutions for \eqref{2-18-4} (or \eqref{3-18-4}). But for the
singularly perturbed problems \eqref{2-18-4} and \eqref{3-18-4}, the
parameter $\varepsilon$ will tend to zero as the number of the
solutions tends to infinity. So, all these results do not give any
multiplicity result for  \eqref{equation} (or \eqref{equation1}).

 In this paper,  we assume that $V(y)$ is radial. That is, $V(y)=V(|y|)$. Thus, we consider the following problem
\begin{equation}
\label{1.4}
 -\Delta u +V(|y|)u = u^p, u >0 \ \  \mbox{in}\  \R^N, \ \ 
  u \in H^1(\R^N),
\end{equation}
where $1<p<\frac{N+2}{N-2}$ if $N\ge 3$, $1<p<+\infty$ if $N=2$. We
assume
\[
\lim_{|y|\to +\infty} V(|y|)=V_0>0.
\]
Note that if $V(r)$ is non-decreasing, by \cite{GNN}, any solution
of \eqref{1.4} is radial.

The aim of this paper is to obtain {\bf infinitely many non-radial
positive solutions} for \eqref{1.4} under an assumption for $V(r)$ near the
infinity. We assume that $V(r)>0$ satisfies
 the following condition:

(V):  There is are constants $a>0$, $m>1$, $\theta>0$, and $V_0>0$,
such that

\begin{equation}\label{V}
V(r)= V_0+\frac a {r^m} +O\bigl(\frac1{r^{m+\theta}}\bigr),
\end{equation}
as $r\to +\infty$.  (Without loss of generality, we may assume that $
V_0=1.$)

Our main result in this paper can be stated as follows:

\medskip
\begin{theorem}
\label{main}  If $V(r)$ satisfies (V), then problem (\ref{1.4}) has
infinitely many non-radial positive solutions.
\end{theorem}

\medskip

\begin{remark}
To obtain the result in Theorem~\ref{main}, \eqref{V} can not be
changed to

\begin{equation}\label{V1}
V(r)= V_0-\frac a {r^m} +O\bigl(\frac1{r^{m+\theta}}\bigr),
\end{equation}
In fact, it is easy to find a function $V(r)$, satisfying $V'(r)\ge
0$ and \eqref{V1}. So, for this $V$, all the solutions must be
radial.

\end{remark}

\begin{remark}\label{r}
The radial symmetry can be replaced by the following  weaker
symmetry assumption: after suitably rotating the coordinate system,

\begin{itemize}
\item[(V1)]
$  V(y)= V(y^{'}, y^{''}) = V(|y^{'}|, |y_{3}|, ..., |y_{N}|), \
\mbox{where} \   y=(y^{'}, y^{''}) \in \R^{2} \times \R^{N-2}$,

\item[(V2)]
$V(y)= V_0+\frac a {|y|^m} +O\bigl(\frac1{|y|^{m+\theta}}\bigr)$
 as $|y|\to +\infty$, where $a>0$, $m>1$,
$\theta>0$, and $V_0>0$ are some constants.

\end{itemize}

\end{remark}

We believe that Theorem~\ref{main} is still true for non-radial
potential $V(y)$.  So,  we make the following conjecture:

\noindent
{\it \textbf{Conjecture:} Problem~\eqref{equation} has infinitely many
solutions, if
\[
V(y)= V_0+\frac a {|y|^m} +O\bigl(\frac1{|y|^{m+\theta}}\bigr),
\]
as $|y|\to +\infty$, where $V_0>0$, $a>0$ , $m>0$ and $\theta>0$ are
some constants.}

\begin{remark}
Using the same argument, we can prove that if
\[
Q(r)= Q_0-\frac a {r^m} +O\bigl(\frac1{r^{m+\theta}}\bigr),
\]
where $Q_0>0$, $ a>0$, $m>1$ and $\theta>0$ are some constants, then
\[
-\Delta u+u=Q(|y|)u^p,\quad u>0,\quad u\in H^1(\R^N),
\]
has infinitely many positive  non-radial solutions.

\end{remark}
\begin{remark}
If $V(r)$ tends to $V_0$  from below, i.e. $V$ satisfies
\begin{equation}\label{V-1}
V(r)= V_0-\frac a {r^m} +O\bigl(\frac1{r^{m+\theta}}\bigr),
\end{equation}
as $r\to +\infty$,  then  we can use similar method to prove the existence of {\em infinitely many sign-changing solutions}. This recovers the result in \cite{CDS}, at least when $V$ is radially symmetric.

\end{remark}

Before we close this introduction, let us outline the main idea in
the proof of Theorem~\ref{main}.

We will construct solutions with large number of bumps near the
infinity. Since we assume
\[
\lim_{|y|\to+\infty} V(|y|)=1,
\]
we will use the solution of

\begin{equation}\label{10-18-4}
\begin{cases}
-\Delta u +u= u^{p},\;\;u>0, &\text{in}\; \R^N,\\
u(y)\to 0, &\text{as}\; |y|\to \infty,
\end{cases}
\end{equation}
to build up the approximate solutions for \eqref{1.4}. It is
well-known   that \eqref{10-18-4} has a unique solution  $U$,
satisfying $U(y)=U(|y|)$, $U'<0$.

Let

\[
x_j=\bigl(r \cos\frac{2(j-1)\pi}k, r\sin\frac{2(j-1)\pi}k,0\bigr),\quad j=1,\cdots,k,
\]
where $0$ is the zero vector in $\R^{N-2}$,  $r\in [r_0 k \ln k, r_1
 k \ln k]$ for some  $r_1>r_0>0$.

Set $y=(y',y'')$, $y'\in \R^2$, $y''\in \R^{N-2}$.  Define
\[
\begin{split}
H_s=\bigl\{ u: & u\in H^1(\R^N), u\;\text{is even in} \;y_h, h=2,\cdots,N,\\
& u(r\cos\theta , r\sin\theta, y'')=
u(r\cos(\theta+\frac{2\pi j}k) , r\sin(\theta+\frac{2\pi j}k), y'')
\bigr\}.
\end{split}
\]

Let
\[
W_r(y)=\sum_{j=1}^k U_{x_j}(y),
\]
where $U_{x_j}(y)=U(y-x_j)$.

Theorem~\ref{main} is a direct consequence of the following result:

\begin{theorem}
\label{th11} Suppose that $V(r)$ satisfies (V). Then there is an
integer $k_0>0$, such that for any integer  $k\ge k_0$, \eqref{1.4}
has a solution $u_k$ of the form
\[
u_k = W_{r_k}(y)+\omega_k,
\]
where  $\omega_k\in H_s$,   $r_k \in [r_0 k \ln k, r_1
 k \ln k]$ and as $k\to +\infty$,
\[
\int_{\R^N} \bigl(|D \omega_k |^2+\omega_k^2\bigr)\to 0.
\]
\end{theorem}

 We will use the techniques in the singularly perturbed elliptic
problems to prove Theorem~\ref{th11}.   For singularly perturbed problems (similar to (\ref{2-18-4})),  a small parameter is present either in the front of the $\Delta$ or in the nonlinearity. Here for our problem (\ref{1.4}), there is no parameter to use. However we use the {\em loss of compactness} to build up  solutions. More precisely, because of the domain $\R^N$ 
we can use $k$, {\bf the number of the
bumps} of the solutions, as the parameter in the construction of
spike solutions for \eqref{1.4}. This seems to be a  {\bf new idea}. This is partly motivated by recent paper of Lin-Ni-Wei 
\cite{LNW} where  they constructed multiple spikes to a singularly perturbed problem.  There they allowed the number of spikes to depend on the small parameter.

This paper is organized as follows. In section~2, we will carry out
the reduction. Then, we will study the reduced finite dimensional
problem and prove Theorem~\ref{th11}. We will leave all the
technical calculations in the appendix.

\medskip

\noindent
{\bf Acknowledgment.} The first
author is supported by an Earmarked
Grant  from RGC  of Hong Kong. The second author is partially supported
by ARC.

\section{Proof the the main result}
\setcounter{equation}{0}

Let
 \[
Z_j=\frac{\partial U_{x_j}}{\partial r},\quad j=1,\cdots,k,
  \]
where $x_j=\bigl(r \cos\frac{2(j-1)\pi}k,
r\sin\frac{2(j-1)\pi}k,0\bigr) $. In this paper, we always assume

\begin{equation}\label{1-20-4}
r\in S_k=: \bigl[ (\frac m{2\pi}-\beta)k\ln k, (\frac m{2\pi}+\beta) k\ln k\bigr],
\end{equation}
where $m$ is the constant in the expansion for $V$, and  $\beta>0$
is a small constant.

Define

\[
E=\bigl\{ v:  v\in H_s,\; \int_{\R^N} U_{x_j}^{p-1} Z_j v=0, \;j=1,\cdots,k\bigr\}
\]
The norm of $H^1(\R^N)$ is defined as follows:

\[
\|v\|=\sqrt{\bigl\langle v, v\bigr\rangle},\quad v\in H^1(\R^N),
\]
where

\[
\bigl\langle v_1, v_2\bigr\rangle=\int_{\R^N} \bigl(Dv_1 Dv_2+V(|y|) v_1 v_2\bigr).
\]

It is easy to check that

\[
 \int_{\R^N} \bigl(Dv_1 Dv_2+V(|y|) v_1 v_2- p W_r^{p-1} v_1 v_2
\bigr),\quad v_1,
v_2 \in E,
\]
is a bounded bilinear functional in $E$. Thus,
  there is a bounded linear operator $L$ from $E$ to $E$, such that

\[
\bigl\langle Lv_1, v_2\bigr\rangle = \int_{\R^N} \bigl(Dv_1 Dv_2+V(|y|) v_1 v_2- p W_r^{p-1} v_1 v_2
\bigr),\quad v_1,
v_2 \in E.
\]

The next lemma shows that $L$ is invertible in $E$.

\begin{lemma}\label{l21}
 There is a constant $\rho>0$, independent of $k$, such that for any
 $r\in S_k$,

 \[
 \|Lv\|\ge \rho \|v\|,\quad  v\in E.
 \]

\end{lemma}

\begin{proof}

We argue by contradiction.  Suppose that there are $n\to +\infty$,
$r_k\in  S_k$, and $v_k\in E$, with
\[ \|Lv_k\|= o(1) \|v_k\|.
 \]
Then

\begin{equation}\label{1-l21}
\bigl\langle L v_k,\varphi\bigr\rangle = o(1) \|v_k\|\|\varphi\|,\quad \forall\; \varphi\in E.
\end{equation}
We may assume that $\|v_k\|^2=k$.

Let

\[
\Omega_j=\bigl\{ y=(y',y'')\in\R^2\times \R^{N-2}:
 \bigl\langle \frac {y'}{|y'|}, \frac{x_j}{|x_j|}\bigr\rangle\ge \cos \frac{\pi}{k}\bigr\}.
\]
By symmetry, we see from \eqref{1-l21},

\begin{equation}\label{2-l21}
\int_{\Omega_1} \bigl(Dv_k D\varphi +V(|y|) v_k \varphi- p W_r^{p-1} v_k \varphi
\bigr)=\frac 1k\bigl\langle L v_k,\varphi\bigr\rangle = o\bigl(\frac1{\sqrt k}\bigr)
\|\varphi\|,\quad \forall\; \varphi\in E.
\end{equation}
In particular,

\[
\int_{\Omega_1} \bigl(|Dv_k|^2+V(|y|) v_k^2- p W_r^{p-1} v_k^2
\bigr)= o(1),
\]
and

\begin{equation}\label{0-19-3}
\int_{\Omega_1} \bigl(|Dv_k|^2+V(|y|) v_k^2\bigr)=1.
\end{equation}

Let  $ \bar v_k(y)= v_k(y-x_1)$. Then for any $R>0$, since
$|x_2-x_1|= r \sin\frac{\pi}k\ge \frac m4 \ln k$, we see that
$B_R(x_1)\subset \Omega_1$.  As a result, from \eqref{0-19-3}, we
find that for any $R>0$,

\[
\int_{B_R(0)} \bigl(|D\bar v_k|^2+V(|y|) \bar v_k^2\bigr)\le 1.
\]
So, we may assume that there is a $v\in H^1(\R^N)$, such that as
$k\to +\infty$,

\[
\bar v_k \to v, \quad \text{weakly in}\; H^1_{loc}(\R^N),
\]
and

\[
\bar v_k \to v, \quad \text{strongly in}\; L^2_{loc}(\R^N).
\]

Since $\bar v_k$  is even in $y_h$, $h=2,\cdots,N$, it is easy to
see that  $v$ is even in $y_h$, $h=2,\cdots,N$. On the other hand,
from

\[
\int_{\R^N}U_{x_1}^{p-1} Z_1 v_k=0,
\]
we obtain

\[
\int_{\R^N} U^{p-1} \frac{\partial U}{\partial x_1} \bar v_k=0.
\]
So, $v$ satisfies

\begin{equation}\label{5-l21}
\int_{\R^N} U^{p-1} \frac{\partial U}{\partial x_1} v=0.
\end{equation}

Now, we claim that $v$ satisfies

\begin{equation}\label{4-l21}
-\Delta v + v -p U^{p-1} v=0,\quad\text{in}\; \R^N.
\end{equation}

Define

\[
\tilde E =\bigl\{ \varphi:  \varphi\in H^1(\R^N),\;\int_{\R^N} U^{p-1}\frac{\partial U}{\partial x_1}\varphi=0
\bigr\}.
\]

For any $R>0$,  let $\varphi\in C_0^\infty(B_R(0))\cap \tilde E$ be
any function, satisfying that  $\varphi$ is even in $y_h$,
$h=2,\cdots,N$.  Then $\varphi_k(y)=: \varphi (y-x_1)\in
C_0^\infty(B_R(x_1))$.  Inserting $\varphi_k$ into \eqref{2-l21},
using Lemma~\ref{al1}, we find

\begin{equation}\label{6-l21}
\int_{\R^N} \bigl( Dv D\varphi +v\varphi -p U^{p-1} v\varphi\bigr)=0.
\end{equation}

On the other hand,  since $v$ is even in $y_h$,  $h=2,\cdots,N$,
\eqref{6-l21} holds for any function $\varphi\in C_0^\infty(\R^N)$,
which is odd  in $y_h$,  $h=2,\cdots,N$. Therefore,  \eqref{6-l21}
holds for any  $\varphi\in C_0^\infty(B_R(0))\cap \tilde E$.  By the
density of $C_0^\infty(\R^N)$ in $H^1(\R^N)$, it is easy to show
that

\begin{equation}\label{7-l21}
\int_{\R^N} \bigl( Dv D\varphi +v\varphi -p U^{p-1} v\varphi\bigr)=0,\quad \forall\; \varphi\in \tilde E.
\end{equation}
  But
\eqref{7-l21} holds for $\varphi =\frac{\partial U}{\partial x_1} $.
Thus \eqref{7-l21} is true for any $\varphi\in H^1(\R^N)$. So, we
have proved \eqref{4-l21}.

Since $U$ is non-degenerate, we see that $v=c \frac{\partial
U}{\partial x_1}$ because $v$ is even in $y_h$, $h=2,\cdots,N$. From
\eqref{5-l21}, we find

\[
v=0.
\]
As a result,

\[
\int_{B_R(x_1)} v_k^2=o(1),\quad \forall\; R>0.
\]

On the other hand, it follows from Lemma~\ref{al1} that for any
small $\eta>0$, there is a constant $C>0$, such that

\begin{equation}\label{1-19-3}
W_{r_k}(y)\le C e^{-(1-\eta)|y-x_1|},\quad y\in\Omega_1.
\end{equation}
Thus,

\[
\begin{split}
o(1)=&\int_{\Omega_1} \bigl(|Dv_k|^2+V(|y|) v_k^2- p W_{r_k}^{p-1} v_k^2
\bigr)\\
=&\int_{\Omega_1} \bigl(|Dv_k|^2+V(|y|) v_k^2\bigr)+ o(1) +O(e^{-(1-\eta)(p-1)R})\int_{\Omega_1}v_k^2\\
\ge &\frac12\int_{\Omega_1} \bigl(|Dv_k|^2+V(|y|) v_k^2\bigr)+ o(1).
\end{split}
\]
This is a contradiction to \eqref{0-19-3}.

\end{proof}

Define

\[
I(u)=\frac12\int_{\R^N} \bigl( |Du|^2+ V(|y|) u^2-\frac1{p+1}\int_{\R^N} |u|^{p+1}.
\]
Let

\[
J(\phi)= I(W_r+\phi),\quad  \phi\in E.
\]

 We have

\begin{proposition}\label{p1-6-3}
There is an integer $k_0>0$, such that for each $k\ge k_0$,   there
is a $C^1$ map from $S_k$ to $H_s$: $\phi=\phi(r)$, $r=|x_1|$,
satisfying $\phi\in E$, and

\[
\bigl\langle \frac{\partial J(\phi)}{\partial \phi},\varphi\bigr\rangle =0,\quad \forall\; \varphi\in E.
\]
Moreover, there is a small $\sigma>0$, such that
\begin{equation}\label{2-20-4}
\|\phi\|\le \frac{C }{k^{\frac {m-1}2+\sigma}}.
\end{equation}

\end{proposition}

\begin{proof}

Expand $J(\phi)$ as follows:

\[
J(\phi)=J(0)+l(\phi)+ \frac12 \bigl\langle L\phi,\phi\bigr\rangle +R(\phi), \quad \phi\in E,
\]
where

\[
l(\phi)=\sum_{i=1}^k \int_{\R^N}(V(|y|)-1) U_{x_i} \phi +\int_{\R^N}
\bigl( W^{p}-\sum_{i=1}^k U_{x_i}^p\bigr)\phi,
\]
$L$ is the bounded linear map from $E$ to $E$ in Lemma~\ref{l21},
and

\[
R(\phi)=\frac1{p+1} \int_{\R^N} \bigl( |W+\phi|^{p+1} -W^{p+1} -(p+1) W^p \phi -\frac12(p+1)p
W^{p-1} \phi^2\bigr).
\]
Since $l(\phi)$ is a bounded linear functional in $E$, we know that
there is an $l_k\in E$, such that

\[
l(\phi)=\bigl\langle l_k,\phi\bigr\rangle.
\]
Thus, finding a critical point for $J(\phi)$ is equivalent to
solving

\begin{equation}\label{1-20-3}
l_k+ L\phi +R'(\phi)=0.
\end{equation}
By Lemma~\ref{l21}, $L$ is invertible.  Thus, \eqref{1-20-3} can be
rewritten as

\[
\phi= A(\phi)=:- L^{-1} l_k - L^{-1} R'(\phi).
\]

Let

\[
S=\bigl\{ \phi:  \phi\in E,  \|\phi\|\le \frac1{k^{\frac {m-1}2} }\bigr\}.
\]

If $p\le 2$, then it is easy to check that

\[
\|R'(\phi)\|\le C\|\phi\|^p.
\]
So, from Lemma~\ref{l1-25-3} below,

\begin{equation}\label{3-20-4}
\|A(\phi)\|\le C\|l_k\|+ C\|\phi\|^p \le
 \frac{C }{k^{\frac {m-1}2+\sigma}}+
 \frac {C}{k^{\frac {p(m-1)}2} }
\le  \frac1{k^{\frac {m-1}2} }.
\end{equation}
Thus, $A$ maps $S$ into $S$ if $p\le 2$.

On the other hand, if $p\le 2$, then

\[
\|R''(\phi)\|\le C\|\phi\|^{p-1}.
\]
Thus,

\[
\begin{split}
&\|A(\phi_1)-A(\phi_2)\| = \| L^{-1} R'(\phi_1)-L^{-1} R'(\phi_2)\|\\
\le & C\bigl( \|\phi_1\|^{p-1}+\|\phi_2\|^{p-1}\bigr)\|\phi_1-\phi_2\|\le \frac12 \|\phi_1-\phi_2\|.
\end{split}
\]
So,  we have proved that if $p\le 2$,   $A$ is a contraction map.
Therefore,  we have proved that if $p\le 2$, $A$ is a contraction
map from $S$ to $S$. So, the result follows from the contraction
mapping theorem.

It remains to deal with the case $p>2$.

Suppose that $p>2$.  Since

\[
\bigl|\bigl\langle R'(\phi),\xi\bigr\rangle\bigr|\le C\int_{\R^N}W_r^{p-2}|\phi|^2|\xi|\le
C\Bigl(\int_{\R^N} \bigl(W_r^{p-2} |\phi|^2\bigr)^{\frac{p+1}p}\Bigr)^{\frac p{p+1}}\|\xi\|.
\]
we find

\[
\|R'(\phi)\|\le C \Bigl(\int_{\R^N} \bigl(W_r^{p-2} |\phi|^2\bigr)^{\frac{p+1}p}\Bigr)^{\frac p{p+1}}.
\]

 On the other hand, it follows from Lemma~\ref{al1} that $W_r$ is
 bounded. Since $2<\frac{2(p+1)}p<p+1$, we obtain

\[
\|R'(\phi)\|
\le C\Bigl(\int_{\R^N}  |\phi|^{\frac{2(p+1)}p}\Bigr)^{\frac p{p+1}}
\le  \|\phi\|^2.
\]

For the estimate of $\|R''(\phi)\|$, we have

\[
\begin{split}
&\bigl|R''(\phi)(\xi,\eta)\bigr|\le C\int_{\R^N} W^{p-2}_r |\phi||\xi||\eta|
\le C\int_{\R^N}  |\phi||\xi||\eta|\\
\le & C \bigl(\int_{\R^N}|\phi|^3 \bigr)^\frac13
 \bigl(\int_{\R^N}|\xi|^3 \bigr)^\frac13  \bigl(\int_{\R^N} |\eta|^3\bigr)^\frac13\le C\|\phi\|\|\xi\|\|\eta\|,
 \end{split}
\]
since $2<3<p+1$. So
\[
\|R''(\phi)\|
\le  C\|\phi\|.
\]
Thus,

\[
\|R'(\phi)\| \le  \frac C{k^{m-1}}\le \frac C {k^{\frac{m-1}2+\sigma}}.
\]
As a result,

\begin{equation}\label{4-20-4}
\|A(\phi)\|\le C\|l_k\|+ C\|R'(\phi)\| \le
 \frac{C }{k^{\frac {m-1}2+\sigma}}
\le  \frac1{k^{\frac {m-1}2} }.
\end{equation}
Thus,  $A$ maps $S$ to $S$.

On the other hand,

\[
\|R''(\phi)\|\le  C\|\phi\| \le \frac {C}{k^{\frac{m-1}2}},
\]
which implies that  $A$ is a contraction map.  So, we have proved
that if $p>2$, then $A$ is a contraction map from $S$ to $S$. And
the result follows from the contraction mapping theorem.

Finally, \eqref{2-20-4} follows from \eqref{3-20-4} and
\eqref{4-20-4}.

\end{proof}

\begin{lemma}\label{l1-25-3}
There is a small $\sigma>0$, such that

\[
\|l_k\|\le \frac{C }{k^{\frac {m-1}2+\sigma}}.
\]

\end{lemma}

\begin{proof}
By the symmetry of the problem,

\begin{equation}\label{1-l10}
\begin{split}
&\sum_{i=1}^k \int_{\R^N}(V(|y|)-1) U_{x_i} \phi=k\int_{\R^N}(V(|y|)-1) U_{x_1} \phi\\
=&k\int_{\R^N}(V(|y-x_1|)-1) U \phi (y-x_1)\le
k O\bigl(\frac1{r^m}\bigr)\|\phi\|\\
\le & \frac{C }{k^{\frac {m-1}2+\sigma}}\|\phi\|,
\end{split}
\end{equation}
because $m>1$.

 On the other hand, for any $y\in \Omega_1$,

 \[
 U_{x_i}^p\le U_{x_1}^{p-1} U_{x_i}.
 \]
Thus,

\begin{equation}\label{2-l10}
\begin{split}
&\Bigl|\int_{\R^N}
\bigl( W^{p}-\sum_{i=1}^k U_{x_i}^p\bigr)\phi\Bigr|=k\Bigl|\int_{\Omega_1}
\bigl( W^{p}-\sum_{i=1}^k U_{x_i}^p\bigr)\phi\Bigr|\\
\le  & C k \int_{\Omega_1} U_{x_1}^{p-1}\sum_{j=2}^k U_{x_j} |\phi|
\le C k  \sum_{j=2}^ke^{-\frac {p-\tau}2|x_j-x_1|}\Bigl(\int_{\Omega_1} |\phi|^{p+1}
\Bigr)^{\frac 1{p+1}}
\\
\le & k^{\frac p{p+1}} \sum_{j=2}^ke^{-\frac {p-\tau}2|x_j-x_1|}\|\phi\|,
\end{split}
\end{equation}
where $\tau>0$ is any small fixed constant.

From  the definition of $S_k$ in  \eqref{1-20-4}, we see that for
any $r\in S_k$

\[
\sum_{j=2}^ke^{-\frac {p-\tau}2|x_j-x_1|}\le Ce^{-\frac{p-\tau}2\frac{2\pi}k r}\le \frac{C}{k^{\frac {p-\tau}2 (m-\beta)
}}.
\]

Since

\[
\frac{pm}2 -\frac{p}{p+1}>\frac{m-1}2,\quad \forall\; m>1,
\]
we obtain from \eqref{2-l10} that

\begin{equation}\label{3-l10}
\Bigl|\int_{\R^N}
\bigl( W^{p}-\sum_{i=1}^k U_{x_i}^p\bigr)\phi\Bigr|
\le  \frac{C }{k^{\frac {m-1}2+\sigma}}\|\phi\|.
\end{equation}
The result follows from \eqref{1-l10} and \eqref{3-l10}.

\end{proof}

We are ready to prove Theorem~\ref{th11}. Let $\phi=\phi(r)$ be the
map obtained in  Proposition~\ref{p1-6-3}. Define

\[
F(r)=I(W_r+\phi),\quad \forall\; r\in S_k.
\]
It is well known that if $r$ is a critical point of $F(r)$, then
$W_r+\phi$ is a solution of \eqref{1.4}. (See \cite{KW} or \cite{LNW}.)

\begin{proof}[Proof of Theorem~\ref{th11}]

 It follows from
Propositions~\ref{p1-6-3} and \ref{p-a1}  that

\[
\begin{split}
&F(r)=I(W)+l(\phi)+ \frac12 \bigl\langle L\phi,\phi\bigr\rangle +R(\phi)\\
=&I(W)+O\bigl(\|l_k\|\|\phi\|+\|\phi\|^2\bigr)\\
=& k\Bigl( A+\frac{B_1}{r^m}- B_2 e^{-\frac {2\pi r}k }+ O\bigl(\frac 1{k^{m+\sigma} }\bigr)
\Bigr)
\end{split}
\]

Consider

\begin{equation}\label{2-25-3}
\max\bigl\{F(r):\; r\in S_k \bigr\}.
\end{equation}
For the definition  of $S_k$, see \eqref{1-20-4}.  Since the
function

\[
\frac{B_1}{r^m}- B_2 e^{-\frac {2\pi r}k }
\]
has a maximum point

\[
\bar r_k= \bigl(\frac m{2\pi}+o(1)\bigr) k \ln k,
\]
which is an interior point of  $S_k$, it is easy to check that
\eqref{2-25-3} is achieved by some $r_k$, which is in the interior
of $S_k$.  Thus, $r_k$ is a critical point of $F(r)$. As a result
\[
W_{r_k}+\phi(r_k)
\]
is a solution of \eqref{1.4}.

\end{proof}

\appendix

\section{Energy Expansion}

In this section, we will give the energy expansion for the
approximate solutions. Recall

\[
x_j=\bigl(r \cos\frac{2(j-1)\pi}k, r\sin\frac{2(j-1)\pi}k,0\bigr),\quad j=1,\cdots,k,
\]

\[
\Omega_j=\bigl\{ y=(y',y'')\in\R^2\times\R^{N-2}:
 \bigl\langle \frac {y'}{|y'|}, \frac{x_j}{|x_j|}\bigr\rangle\ge \cos\frac\pi k\bigr\},
\quad j=1,\cdots,k,
\]
and

\[
I(u)=\frac12\int_{\R^N}\bigl(|Du|^2+u^2\bigr)-\frac1{p+1}\int_{\R^N} |u|^{p+1}.
\]

Firstly, we have the following basic estimate:

\begin{lemma}\label{al1}
For any $y\in\Omega_1$, and $\eta\in (0,1]$, there is a constant
$C>0$, such that

\[
\sum_{j=2}^k U_{x_j}(y)\le C e^{-\eta |x_1|\frac{\pi}k}e^{-(1-\eta)|y-x_1|}.
\]

\end{lemma}

\begin{proof}

For any    $y\in\Omega_1$, we have
  $|y-x_j|\ge |y-x_1|$.

If $|y-x_1|\ge 2|x_j-x_1|$, then for  any    $y\in\Omega_1$,

\[
\begin{split}
&U_{x_j}(y)\le Ce^{-|y-x_j|}\le Ce^{-|y-x_1|}\\
=& Ce^{-\eta|y-x_1|}e^{-(1-\eta)|y-x_1|}\le  Ce^{-2\eta|x_j-x_1|}e^{-(1-\eta)|y-x_1|}.
\end{split}
\]

If $|y-x_1|\le 2|x_j-x_1|$, then

\[
|y-x_j|\ge |x_j-x_1|-|y-x_1|\ge \frac12|x_j-x_1|.
\]
So for  any    $y\in\Omega_1$,

\[
U_{x_j}(y)\le C e^{-\eta|y-x_j|}e^{-(1-\eta)|y-x_1|}\le Ce^{-\frac12\eta|x_j-x_1|}e^{-(1-\eta)|y-x_1|}.
\]
Thus,

\[
\begin{split}
&\sum_{j=2}^k U_{x_j}(y)\le C e^{-(1-\eta)|y-x_1|} \sum_{j=2}^k e^{-\frac12\eta|x_j-x_1|}\\
\le & C e^{-(1-\eta)|y-x_1|} \sum_{j=2}^k e^{-\eta |x_1|\sin\frac{j\pi}k}
\le C_1 e^{-\eta |x_1|\frac{\pi}k}e^{-(1-\eta)|y-x_1|}.
\end{split}
\]

\end{proof}

In this appendix, we denote  $r=|x_1|$, and we always assume that

\[
r\in S_k,
\]
where $S_k$ is defined in \eqref{1-20-4}.

\begin{proposition}\label{p-a0}

We have

\[
I(U_{x_1})=A+\frac{B_1}{r^m}+ O\bigl(\frac 1{k^{m+\theta} }\bigr),
\]
where

\[
A=\bigl(\frac12-\frac1{p+1}\bigr) \int_{\R^N} U^{p+1},
\quad  B_1=\frac a2\int_{\R^N} U^2.
\]

\end{proposition}

\begin{proof}

We have

\begin{equation}\label{1-p-a0}
\begin{split}
&I(U_{x_1})=\frac12 \int_{\R^N} \bigl(|D U|^2 +U^2\bigr)-\frac1{p+1}\int_{\R^N} U^{p+1}+\frac12
\int_{\R^N} \bigl( V(|y|)-1\bigr) U^2_{x_1}\\
=& A +\frac12 \int_{\R^N} \bigl( V(|y-x_1|)-1\bigr) U^2.
\end{split}
\end{equation}

On the other hand, for any small $\tau>0$, using (V),

\begin{equation}\label{2-p-a0}
\begin{split}
&\int_{\R^N} \bigl( V(|y-x_1|)-1\bigr) U^2=\int_{B_{\frac12 r}(0)} \bigl( V(|y-x_1|)-1\bigr) U^2+
O\bigl( e^{-(1-\tau) r}\bigr)\\
=&\int_{B_{\frac12 r}(0)} \Bigl( \frac a{|y-x_1|^m}+ O\bigl( \frac 1{|y-x_1|^{m+\theta}}\bigr)\Bigr) U^2+
O\bigl( e^{-(1-\tau) r}\bigr).
\end{split}
\end{equation}
But for any $\alpha>0$,

\[
\frac 1{|y-x_1|^\alpha}=\frac 1{|x_1|^\alpha}\Bigl( 1+O\bigl(\frac{|y|}{|x_1|}\bigr)\Bigr),\quad
y\in B_{\frac12 |x_1|}(0).
\]
Thus,

\begin{equation}\label{3-p-a0}
\int_{B_{\frac12 r}(0)} \frac 1{|y-x_1|^\alpha} U^2= \frac 1{|x_1|^\alpha}\int_{\R^N} U^2+
O\bigl(\frac 1{|x_1|^{\alpha+1}}+ e^{-(1-\tau)|x_1|}\bigr).
\end{equation}

Inserting \eqref{3-p-a0} into \eqref{2-p-a0}, we obtain

\begin{equation}\label{4-p-a0}
\int_{\R^N} \bigl( V(|y-x_1|)-1\bigr) U^2=
\frac{B_1}{r^m}+ O\bigl(\frac 1{k^{m+\theta} \ln^{m+\theta}k}\bigr).
\end{equation}
Thus, the result follows from \eqref{1-p-a0} and \eqref{4-p-a0}.

\end{proof}

\begin{proposition}\label{p-a1}

There is a small constant $\sigma>0$, such that

\[
I(W_r)=k\Bigl( A+\frac{B_1}{r^m}- B_2 e^{-\frac {2\pi r}k }+ O\bigl(\frac 1{k^{m+\sigma} }\bigr)
\Bigr),
\]
where  $A$ and $ B_1$ are the constants in Proposition~\ref{p-a0},
and $B_2>0$ is a  positive constant.

\end{proposition}

\begin{proof}

Using the symmetry,

\begin{equation}\label{1-p-a1}
\begin{split}
&\int_{\R^N} (|DW_r|^2 +W_r^2) =  \sum_{j=1}^k \sum_{i=1}^k \int_{\R^N} U_{x_j}^p U_{x_i}\\
=& k \int_{\R^N} U^{p+1}+ k\sum_{i=2}^k \int_{\R^N}U_{x_1}^p U_{x_i}.
\end{split}
\end{equation}

Recall

\[
\Omega_j=\bigl\{ y=(y',y'')\in \R^2\times \R^{N-2}:
 \bigl\langle \frac {y'}{|y'|}, \frac{x_j}{|x_j|}\bigr\rangle \ge \cos\frac\pi k\bigr\},
\quad j=1,\cdots,k.
\]

It follows from Lemma~\ref{al1} that

\begin{equation}\label{2-p-a1}
\begin{split}
&\int_{\R^N} (V(|y|)-1)W_r^2 = k\int_{\Omega_1} (V(|y|)-1)W_r^2 \\
=& k \int_{\Omega_1} (V(|y|)-1)\Bigl( U_{x_1}+O\bigl(e^{-\frac12 |x_1|\frac \pi k}e^{-\frac12 |y-x_1|}\bigr)\Bigr)^2\\
=& k\int_{\Omega_1} (V(|y|)-1)U_{x_1}^2+kO\Bigl(\int_{\Omega_1} |V(|y|)-1|
e^{- \frac12 |x_1|\frac \pi k}e^{- |y-x_1|}
\Bigr)\\
=&
k\Bigl( \frac{B_1}{|x_1|^m} +O\bigl(\frac1{k^{m+\theta}}\bigr)\Bigr).
\end{split}
\end{equation}

Suppose that  $p\le 3$. Then, for any $y\in\Omega_1$,

\[
W_r^{p+1}= U_{x_1}^{p+1} +(p+1) U_{x_1}^{p} \sum_{j=2}^k U_{x_j}+O\Bigl(U_{x_1}^{\frac{p+1}2}
\bigl(\sum_{j=2}^k U_{x_j}\bigr)^{\frac{p+1}2}\Bigr).
\]
Using Lemma~\ref{al1}, we have

\[
\begin{split}
&U_{x_1}^{\frac{p+1}2}
\bigl(\sum_{j=2}^k U_{x_j}\bigr)^{\frac{p+1}2}=U_{x_1}^{\frac{p+1}2}\sum_{j=2}^k U_{x_j}
\bigl(\sum_{j=2}^k U_{x_j}\bigr)^{\frac{p-1}2}\\
\le & C e^{-\eta\frac{p-1}2\frac{|x_1|\pi}k}
U_{x_1}^{(1-\eta)p}\sum_{j=2}^k U_{x_j}
\end{split}
\]

But for any $r\in S_k$,

\[
\sum_{j=2}^k e^{-|x_j-x_1|}\le Ce^{-\frac{2\pi}k r}\le \frac C{k^{m-\beta}}.
\]
So, we obtain that for $p\in (1,3]$, if $\beta>0$ is small enough,

\begin{equation}\label{3-p-a1}
\begin{split}
&\int_{\R^N} W_r^{p+1} = k\int_{\Omega_1} W_r^{p+1} \\
=& k \int_{\Omega_1} \Bigl( U^{p+1}_{x_1}+(p+1)\sum_{i=2}^k U_{x_1}^p U_{x_i}\Bigr) +
 O\bigl(\sum_{i=2}^k e^{-\eta\frac{p-1}2\frac{|x_1|\pi}k}e^{-|x_i-x_1|}\bigr)
\Bigr)\\
=&k \Bigl(\int_{\R^N} U_{x_1}^{p+1}+ (p+1)\sum_{i=2}^k \int_{\R^N}U_{x_1}^p U_{x_i}+
 O\bigl(\frac1{k^{m+\sigma}}\bigr)
\Bigr).
\end{split}
\end{equation}

Suppose that $p>3$. Then for any $y\in\Omega_1$,

\[
W_r^{p+1}= U_{x_1}^{p+1} +(p+1) U_{x_1}^{p} \sum_{j=2}^k U_{x_j}+O\Bigl(U_{x_1}^{p-1}
\bigl(\sum_{j=2}^k U_{x_j}\bigr)^{2}\Bigr).
\]
Since $p-1>2$, similar to the proof of \eqref{3-p-a1}, we can obtain
the following estimate for $p>3$:

\begin{equation}\label{4-p-a1}
\begin{split}
&\int_{\R^N} W_r^{p+1} = k\int_{\Omega_1} W_r^{p+1} \\
=&k \Bigl(\int_{\R^N} U_{x_1}^{p+1}+(p+1) \sum_{i=2}^k \int_{\R^N}U_{x_1}^p U_{x_i}+
 O\bigl(\frac1{k^{m+\sigma}}\bigr)
\Bigr).
\end{split}
\end{equation}

Combining \eqref{1-p-a1}, \eqref{2-p-a1}, \eqref{3-p-a1} and
\eqref{4-p-a1}, we are led to

\[
I(W_r)=k\Bigl( A+\frac{B_1}{r^m}-\frac12\sum_{i=2}^k \int_{\R^N}U_{x_1}^p U_{x_i}+
O\bigl(\frac1{k^{m+\sigma}}\bigr)
\Bigr).
\]
But there are constants $\sigma>0$, $\tilde B_2>0$ and $B'_2>0$,
such that

\[
\begin{split}
&\sum_{i=2}^k \int_{\R^N}U_{x_1}^p U_{x_i}= \tilde B_2 \sum_{i=2}^k e^{-|x_i-x_1|} +
O\bigl(\sum_{i=2}^k e^{-(1+\sigma)|x_i-x_1|}\bigr)\\
= &B_2' e^{-\frac{2|x_1|}k} +O\bigl(e^{-\frac{2(1+\sigma)|x_1|}k}\bigr).
\end{split}
\]
So the result follows.

\end{proof}

\end{document}